\documentclass[12pt,letterpaper]{article}
\usepackage[a4paper, total={7in, 10in}]{geometry}

\usepackage{graphicx}
\usepackage{helvet}
\usepackage{authblk}
\usepackage{subcaption}
\usepackage{hyperref}
\usepackage{amsmath} 
\usepackage{amssymb} 
\usepackage{cleveref}
\usepackage{orcidlink} 
\usepackage[super,comma,sort&compress]  
   {natbib}
\usepackage[right]{lineno} \linenumbers

\newtheorem{definition}{Definition}

\newcommand*\chem[1]{\ensuremath{\mathrm{#1}}}

\makeatletter
\renewcommand{\maketitle}{\bgroup\setlength{\parindent}{0pt}
\begin{flushleft}
  \textbf{\@title}
  
  \@author
\end{flushleft}\egroup}
\makeatother

\title{Locational Marginal Emissions for Carbon-Aware Data Center Operations in Large-Scale Power Grids}
\date{}

\author[1]{Luc Cote}
\author[1, 2*]{Andy Sun}

\affil[1]{Operations Research Center, Massachusetts Institute of Technology, Cambridge, MA, USA}
\affil[2]{Sloan School of Management, Massachusetts Institute of Technology, Cambridge, MA, USA}

\affil[*]{Correspondence: sunx@mit.edu}

\begin{document}

\maketitle

\section*{SUMMARY}



Carbon accounting methods for electricity consumption face challenges regarding physical deliverability, double counting, additionality, and impact magnitude. Locational Marginal Emissions (LMEs) show potential to address many of these key issues. However, their use in a large-scale power grids remains understudied. We analyze the properties of LMEs from a data center's perspective in a 1493-bus Western Interconnection over one year of hourly operation. We find that LME characteristics create three distinct regions: the hydropower-dominated Pacific Northwest, with low and stable LMEs; the coal-heavy Intermountain West, containing often high LMEs; and the Sunbelt, where mixed generation leads to variable LMEs correlated with solar output. This characterization provides analytical guidance for data center emission reduction. In particular, LME-guided emission reduction interventions through data center temporal-spatial load shifting, siting, and renewable procurement display over 85\% accuracy with respect to actual emission reduction. Moreover, large-scale, nodal grid simulation is shown to be critical to accurate evaluation.

\section*{KEYWORDS}


Locational marginal emission, carbon accounting, large-scale power systems, data centers

\section*{INTRODUCTION}

As the modern AI economy unfolds, worldwide energy grids are experiencing a surge in data center energy demand just as many nations have begun to reconfigure their energy grids for the renewable transition. Without proper alignment of this massive new data center demand and green grid development, countries risk missing key emissions targets and undoing global climate progress. Indeed, in recent years, countries experiencing data center growth have seen an increase in the energy demand of these facilities outpace the addition of renewable capacity to their grids, contributing to recent failures to reach EU climate targets \cite{cornish_inside_2025, epa_irelands_2025}. 

In spite of such misalignment, many of the same companies driving this increase in demand have also shown a commitment to reduced-emission goals in the form of net-zero targets. To achieve these net-zero goals, data center operators can rely on three major forms of intervention: 
\begin{enumerate}
    \item \textbf{data center siting} - the intentional placement of data centers in areas where low-carbon energy resources are readily available to meet demand,
    \item \textbf{renewable procurement} - the funding of grid upgrade and low carbon generation projects to displace carbon emissions, and
    \item \textbf{data center operations} - the geographic and temporal shift of data center workloads toward areas and periods of low carbon production.
\end{enumerate}

Data center siting and renewable procurement are both types of facility siting decisions and have seen implementation through different forms of scope 2\footnote{scope 2 emissions refer to the attribution of indirect emissions to the consumption of electricity, heating, and cooling.} carbon accounting practices, with data center companies adopting multiple types of carbon accounting policies. In fact, many of the largest data center stakeholders, including Amazon\cite{amazon_carbon-free_2025}, Google\cite{google_clean_2016}, Meta\cite{meta_data_2024}, and Microsoft\cite{microsoft_measuring_2023} publicize their renewable energy procurement and data center location strategies. Furthermore, recent demand-shifting strategies involving data center-driven work loads have shown an ability to participate significantly in demand response programs both in academic studies\cite{cupelli_data_2018} and in practice\cite{giacobone_google_2025}.

However, although such scope 2 interventions enable companies to achieve net zero carbon goals, it is not clear that such interventions reliably work to decrease carbon emissions. Current scope 2 emissions accounting practices have seen significant criticism in recent years, stemming from challenges regarding deliverability, double counting, additionality, and impact magnitude \cite{gillenwater_redefining_2008, brander_creative_2018, bjorn_renewable_2022, Brander, xu_system-level_2024}. Furthermore, recent work has shown data center demand response can lead to negligible emission changes or even increases in system emissions \cite{gorka_electricityemissionsjl_2025, knittel_flexible_2025}.

\paragraph{Locational Marginal Emissions}
One tool that has been identified to provide possible improvements in intervention choice is the use of locational marginal emissions (hereafter, LMEs) which describe the impact of changes in load on the system-wide carbon emissions within an electric power grid\cite{ghgscope2}. More formally, we define the LME corresponding to node $n$ in a power grid as
\begin{definition}
    \begin{equation}
        \text{LME}_n := \frac{d \text{ System Emissions}}{d \text{ Load}_n},
    \end{equation}
\end{definition}
where $\text{System Emissions}$ refers to the \chem{CO_2} output of the entire energy grid and $\text{ Load}_n$ refers to the net power load at node $n$, which is the difference of load and generation at that node. The physical unit of $\text{LME}_n$ is MMT \chem{CO_2}/MWh or kg\chem{CO_2}/MWh \footnote{MMT stands for million metric tons, kg stands for kilograms, and MWh stands for megawatt hours.}. LMEs have also been used to define carbon accounting systems that tie together emission attribution and intervention evaluation \cite{rudkevich_locational_2012,avila_catalyzing_2025} with the following carbon accounting scheme:
\begin{definition} \label{def:accts}
Carbon accounts, 
    \begin{itemize}
    \item \textbf{load:} $\text{LME}_n \times \text{load}_n$ - the carbon account for a load is defined as the nodal marginal emissions rate multiplied by the load's total magnitude,
    \item \textbf{generator:} $(\text{carbon intensity}_g - \text{LME}_g) \times \text{generation}_g$  - the carbon account for a generation unit is defined as the difference between the generator's unit emission rate and the nodal marginal emissions rate multiplied by the total generation,
    \item \textbf{transmission:} $SCI_\ell \times \text{line flow}_{\ell}$ - here, $SCI_\ell$ refers to the ``shadow carbon intensity'' of a transmission line which describes how the system emissions would change with respect to a change in the line's power limit and $\text{line flow}_{\ell}$ refers to the magnitude of the power delivered through line $\ell$.
\end{itemize}
\end{definition}
%
LME-based accounting methods have been shown to have many desirable properties regarding the issues present in current scope 2 accounting practices. Of particular note is the aptly named ``Carbon Accounting Theorem''\cite{rudkevich_locational_2012,avila_catalyzing_2025} which states that the total emissions accounts across all load, generation, and transmission entities as described in \Cref{def:accts} sum to the total scope 1 emissions of the energy grid, allowing such accounting schemes to address double counting issues.
\begin{definition}
    Carbon accounting theorem,
    \begin{align}
        \sum_n \text{LME}_n \times \text{load}_n + \sum_g (\text{carbon intensity}_g - \text{LME}_g) \times \text{generation}_g + \sum_{\ell} \text{SCI}_\ell \times \text{line flow}_{\ell} \notag \\ 
        = \text{system emissions}.
    \end{align}
\end{definition}
Furthermore, LMEs are nodal properties whose definition is based on the response of the grid, and subsequently, the grid congestion that constrains this response is closely related to this calculation \cite{rudkevich_locational_2012,avila_catalyzing_2025}, enabling LME-based accounting to resolve deliverability issues.

\paragraph{Contributions} In this paper, we seek to evaluate the LME characteristics of a large interconnection in the United States and answer how a data center operator might use LMEs and LME-based accounting to guide emission reduction interventions. In particular, we make the following two main contributions:
\begin{enumerate}
    \item investigate broad carbon intervention insights given by the empirical properties of LMEs and LME-based accounting using a high-granularity representation of the WECC grid,
    \item evaluate the accuracy and utility of LMEs as indicators within the WECC grid for the resulting emission changes from interventions of data center siting, renewable procurement, and data center operations.
\end{enumerate}

\textit{Empirical Properties:} We apply a novel optimization-based LME calculation methodology to a 1493-bus realistic representation of the WECC grid with hourly dispatches simulated over 1 year of generation and load data, providing a fine-grained, interconnection-scale survey of grid LMEs. Through such experiments, we analyze the regional, hourly, and seasonal patterns of WECC LMEs and how these patterns are reflected across LME-based carbon accounts. Furthermore, we consider how such results provide insights for a carbon-conscious energy consumer's operations.

\textit{Accuracy Evaluation:}
We extend our simulation experiments to consider the three main LME-based intervention strategies for data centers of consumption siting, renewable procurement, and data center operations. For operations and consumption siting, we consider the lens of a data center operator with large data centers of 100 and 200 MW spanning the WECC grid, while for renewable procurement, we take the perspective of a renewable energy developer performing generation siting of a 100MW zero-carbon resource. In all three cases, we find that LMEs serve as an accurate metric for describing intervention results.

\section*{RESULTS}


In the following experiments, we use a realistic representation of the modern WECC grid to investigate the empirical properties of the described LME calculation and accounting methodology. To create such a representation of the WECC grid, we follow the methodology outlined in \citet{thomas_paper}, using Open Street Map data\cite{OpenStreetMap} along with 2023 EIA 860\cite{EIA860} generation and substation location data to estimate the transmission properties of the WECC 230 KV grid. This leads to a grid with the characteristics described in \Cref{tab:wecc_data}. Additionally, we use historic load and weather data, once again defined in the manner of \citet{thomas_paper}, to define a representative year of WECC hourly dispatch periods, each period defined by scaling factors for nodal loads and variable renewable energy (VRE) capacity factors. The weather data describes the hourly conditions of 2016 at 20km resolution and is sourced from NOAA \cite{noaa_rap252}, while the load data is sourced from the 2023 EIA 930 dataset \cite{EIA930}, which describes the hourly profiles across the WECC area.

\begin{table}
    \centering
    \begin{tabular}{ |c||c|c|c|c|c|c|c|c|c| } 
     \hline
      Entity & Bus & Line & Load & Coal & Gas & Wind & Solar & Nuclear & Hydro \\ 
        \hline
      Number of elements & 1493 & 1919 & 722 & 25 & 180 & 90 & 58 & 3 & 70\\ 
      Installed Capacity (GW) & N/A & 4357.1 & 122.5 & 22.5 & 86.6 & 25.9 & 21.5 & 7.7 & 37.9\\ 
      Emissions ( kg\chem{CO_2}/MWh) & N/A & N/A & N/A & 800 & 450 & 0 & 0 & 0 & 0 \\
     \hline
    \end{tabular}
    \caption{Makeup of the simulated WECC grid}
    \label{tab:wecc_data}
\end{table}

\begin{figure}[h!]
    \centering
    \includegraphics[width=0.75\linewidth]{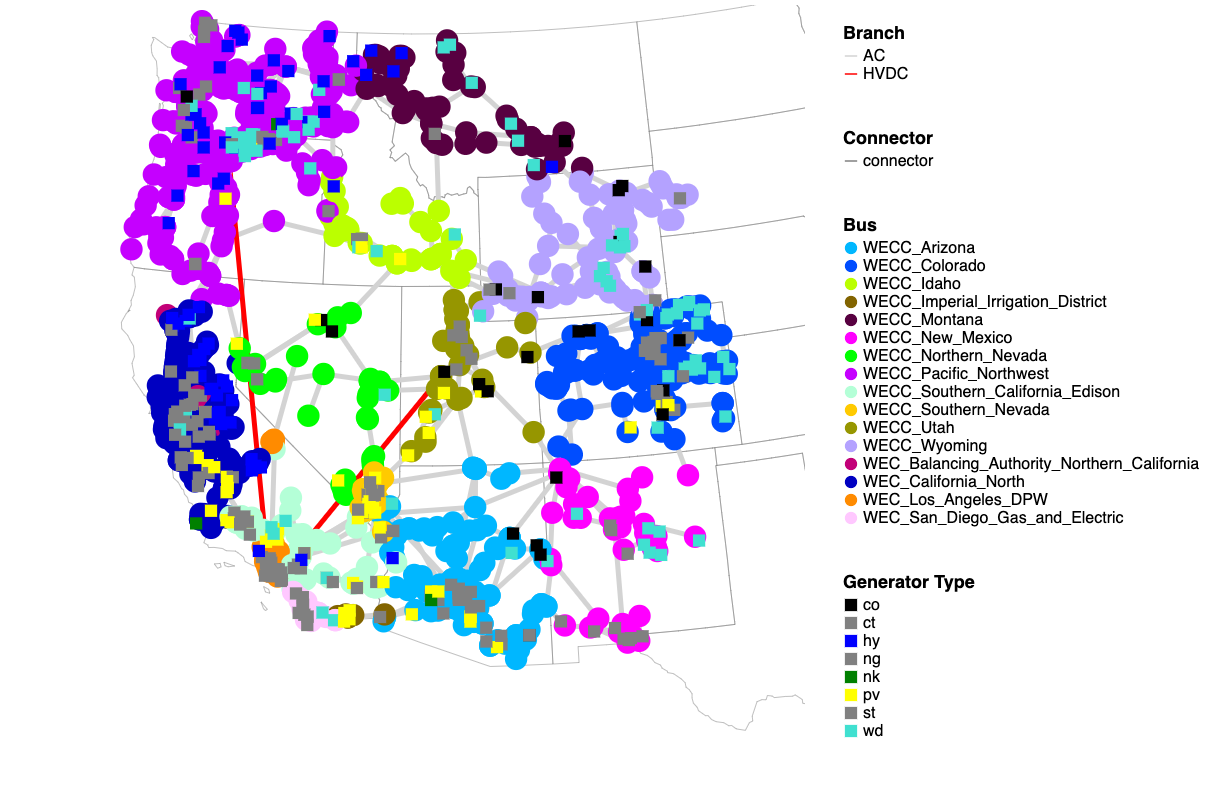}
    \caption{EPA IPM Planning Regions over the WECC Data}
    \label{fig:wecc_regions}
\end{figure}

To observe and quantify regional trends, we additionally refer to the EPA IPM regions across the WECC. In \Cref{fig:wecc_regions}, we visualize how each region corresponds to the buses within the simulated WECC grid.

\subsection*{LMEs Reveal Three Regional Characteristics within WECC} \label{sec:spatial}
\begin{figure}[h!]
    \centering
    \includegraphics[width=0.6\linewidth]{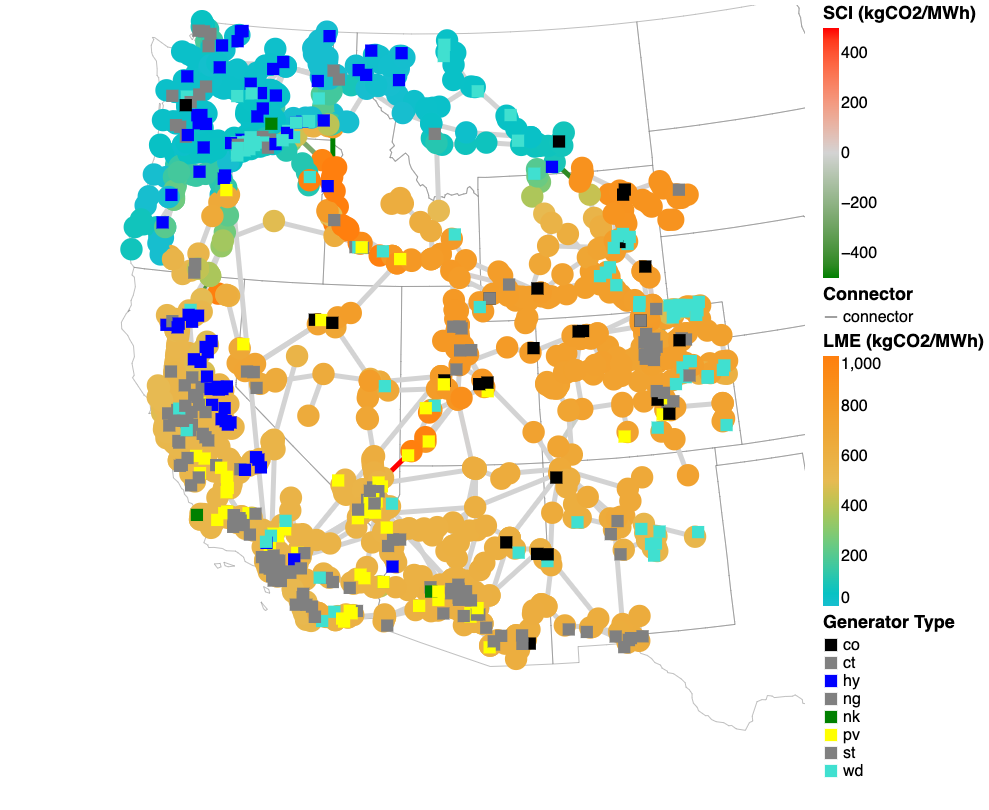}
    \caption{Average locational marginal emissions and shadow carbon intensities over the experimental horizon.}
    \label{fig:wecc_avg_LME}
\end{figure}

\begin{figure}[h!]
    \centering
    \includegraphics[width=0.6\linewidth]{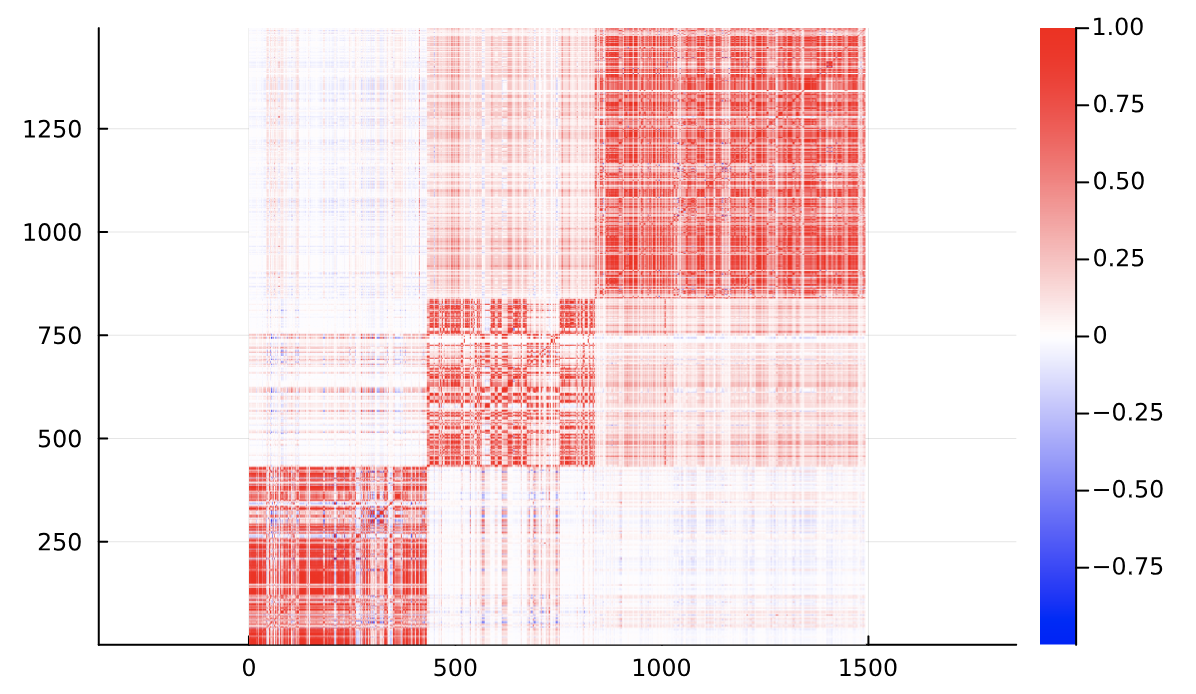}
    \caption{Heatmap of the Pearson correlation coefficients between nodal LMEs. Nodes 1-431 belong to the ``Mountain West'' regions, nodes 432-837 belong to the Pacific Northwest and Montana regions, and nodes 838-1493 belong to ``Sunbelt'' regions.}
    \label{fig:wecc_correlation}
\end{figure}

The WECC grid spans a wide geographic area with a diverse set of climates and policies, leading to significant differences in energy demands and generation resources across regions. Broadly speaking, our grid simulation reveals that WECC regions can be grouped into three main types of LME characteristics whose correlation strength we evaluate in \Cref{fig:wecc_correlation}. 

Pacific Northwest and Montana IPM Regions make up the blue sections in the north of \Cref{fig:wecc_avg_LME}, experiencing consistently low LMEs with an overall area average LME of $\approx$90 kg\chem{CO_2}/MWh. This low LME pattern comes from a high proportion of time with zero carbon generation acting as the marginal generation resource due to a stable supply of hydropower along with significant amounts of wind capacity. An exception to this trend occurs in the southern portion of region, where a less abundant hydro supply combined with some transmission congestion from the north causes nearby natural gas plants to act as the marginal resource at times.

The ``Inter-Mountain West'' IPM regions of Idaho, Wyoming, Colorado, and Utah make up the clusters of darker orange nodes in the center and east of \Cref{fig:wecc_avg_LME}, exhibiting high average LMEs with an area overall average of $\approx$710 kg\chem{CO_2}/MWh. These regions contain some renewable development in the forms of wind and solar; however, they lack sufficient capacity to completely satisfy the regional load during most hours, leading to a tendency for fossil fuel generation to act as the marginal resource. Furthermore, these regions contain significant coal plant capacity, which exhibits emission rates much higher than those of natural gas plants, resulting in high regional LMEs. This region experiences some levels of transmission congestion within the Cheyenne-Denver corridor, where transmission lines between nearby wind farms and the urban areas exhibit negative SCIs\footnote{SCIs refers to the change in system carbon emissions resulting from a change in transmission line capacity as described in \Cref{def:accts}.}, implying congestion limiting wind dispatch. Similarly, the lines connecting the hydropower and wind resources of the PNW and Montana to this region tend to exhibit negative LMEs, reflected in the dark green lines connecting western Washington to Idaho and Western Montana to Wyoming in \Cref{fig:wecc_avg_LME}.

Finally, the remaining ``Sunbelt'' IPM regions make up the lighter orange nodes in the south and west of \Cref{fig:wecc_avg_LME}, exhibiting variable LMEs which follow trends related to solar radiation strength. Such trends can reasonably be explained from high levels of solar production within the regions. Some areas within these regions, such as northern and southern California, contained limited access to coal plants and higher levels of solar development, leading to lower average LMEs along with more pronounced solar trends, compared to other areas such as New Mexico which contain higher coal access along with a more limited solar supply. These regions exhibited mixed congestion with some lines connecting renewables to high demand areas such as from the low-carbon PNW to northern California experiencing negative SCIs, but also lines connecting to coal generation in Utah and New Mexico exhibiting positive SCIs (see the red line northeast of southern Nevada in \Cref{fig:wecc_avg_LME}), indicating the existence of some coal-limiting congestion leading to an increase in natural gas dispatch.

\begin{figure}[h!]
   \begin{subfigure}{0.48\textwidth}
     \centering
     \includegraphics[width=\linewidth]{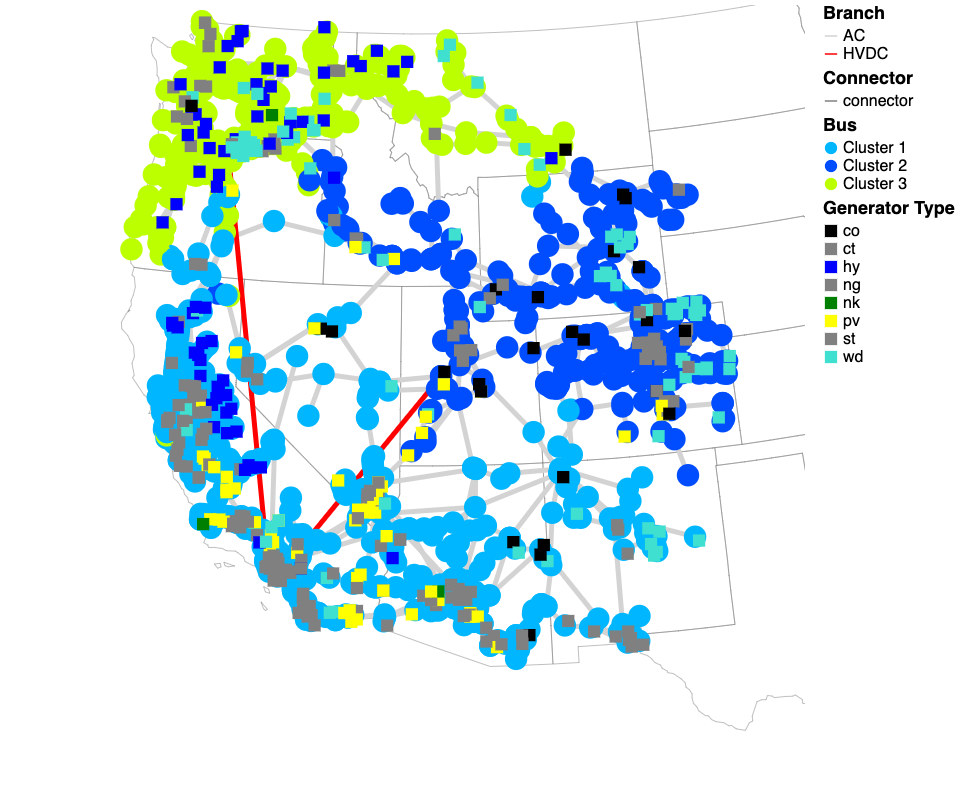}
     \caption{Bus groups recovered for $k=3$ groups.}
    \label{fig:wecc_clustering}
   \end{subfigure}\hfill
   \begin{subfigure}{0.48\textwidth}
     \centering
    \includegraphics[width=\linewidth]{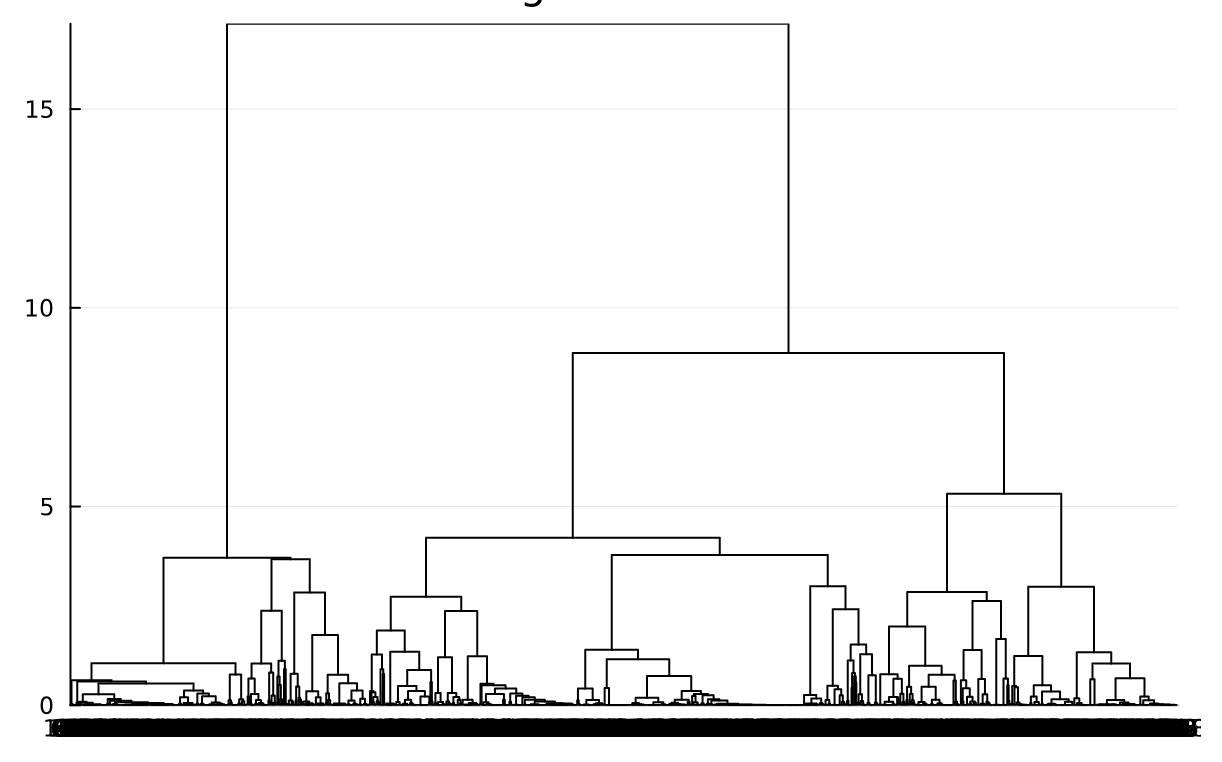}
    \caption{Clustering dendogram describing total internal variance across all clusters.}
    \label{fig:cluster_dendogram}
   \end{subfigure}
   \caption{Results from hierarchical clustering of bus groups using distance \eqref{eq:corr_dist} and ward's method for linkage}
\end{figure}

This three-region grouping of IPM regions is further reflected through observing the correlations between nodal vectors of LMEs over time. \Cref{fig:wecc_correlation} displays a clear three block diagonal structure outlining high intra-regional correlation. Furthermore, we can use such correlation data to create a distance function between buses in the form of 
\begin{equation}
    d_{ij}=1-|\text{correlation}(i,j)|,
\end{equation} to give a measure of the temporal pattern distance between nodes. Combining this distance metric with the rescaled difference magnitude between bus average LMEs allows for the creation of a distance function which reflects the LME characteristics over time 
\begin{equation} 
    \tilde{d}_{ij}=\frac{1}{2}d_{ij}+\frac{1}{2}\frac{|\bar{\text{LME}}_{\text{i}}-\bar{\text{LME}}_{\text{j}}|}{\max_{i',j'}|\bar{\text{LME}}_{\text{i'}}-\bar{\text{LME}}_{\text{j'}}|}, \label{eq:corr_dist}
\end{equation}
which can be used to perform hierarchical clustering among buses \cite{benevento_correlation-based_2024}. Indeed, in \Cref{fig:wecc_clustering} such a hierarchical clustering recovers a bus grouping quite similar to the three region IPM grouping, while in \Cref{fig:cluster_dendogram} the dendrogram created from such a clustering confirms $k=2$ and $k=3$ groups to be natural clustering choices.

\subsection*{LMEs Reflect Hourly and Seasonal Generation Characteristics}\label{sec:temporal}

\begin{figure}[h!]
   \begin{subfigure}{0.48\textwidth}
     \centering
     \includegraphics[width=\linewidth]{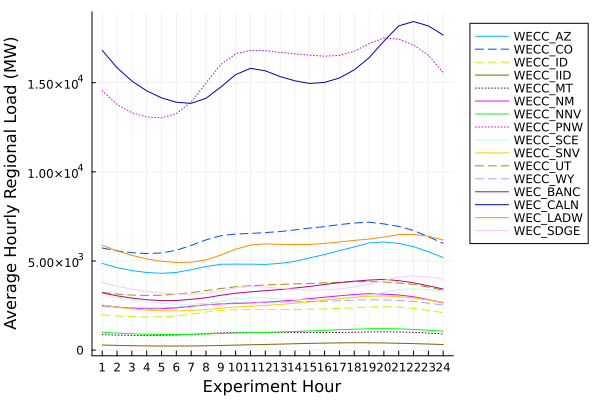}
     \caption{hourly demand}
    \label{fig:hourly_demand}
   \end{subfigure}\hfill
   \begin{subfigure}{0.48\textwidth}
     \centering
    \includegraphics[width=\linewidth]{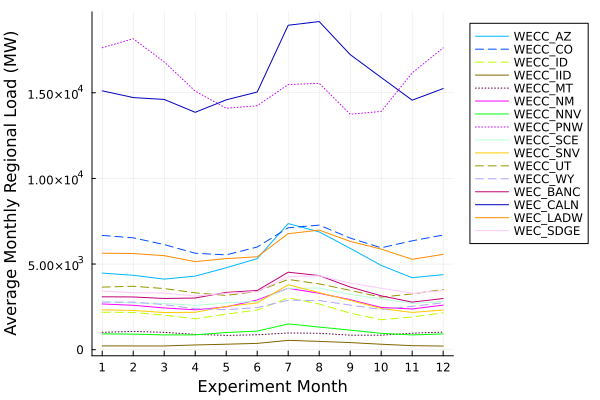}
    \caption{monthly demand}
    \label{fig:monthly_demand}
   \end{subfigure}
   \caption{Average power demand of the WECC EPA IPM regions over hourly and monthly periods.}
\end{figure}

\begin{figure}[h!]
   \begin{subfigure}{0.48\textwidth}
     \centering
     \includegraphics[width=\linewidth]{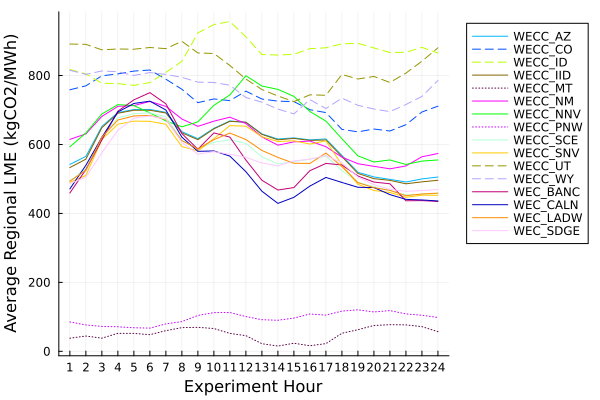}
     \caption{hourly LMEs}
     \label{fig:hourly_LME}
   \end{subfigure}\hfill
   \begin{subfigure}{0.48\textwidth}
     \centering
    \includegraphics[width=\linewidth]{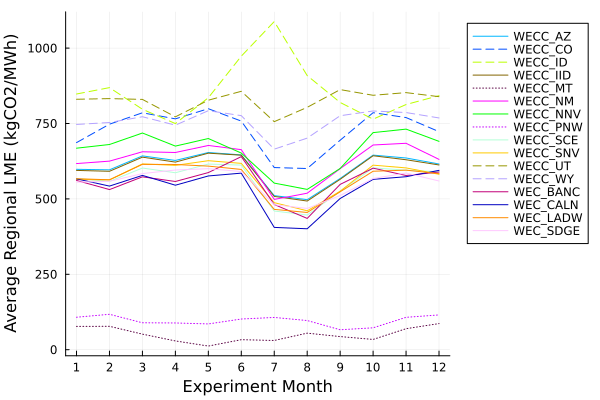}
    \caption{monthly LMEs}
    \label{fig:monthly_LME}
   \end{subfigure}
   \caption{Average locational marginal emissions of the WECC EPA IPM regions over hourly and monthly periods.}
\end{figure}

\begin{figure}[h!]
   \begin{minipage}{0.48\textwidth}
     \centering
     \includegraphics[width=\linewidth]{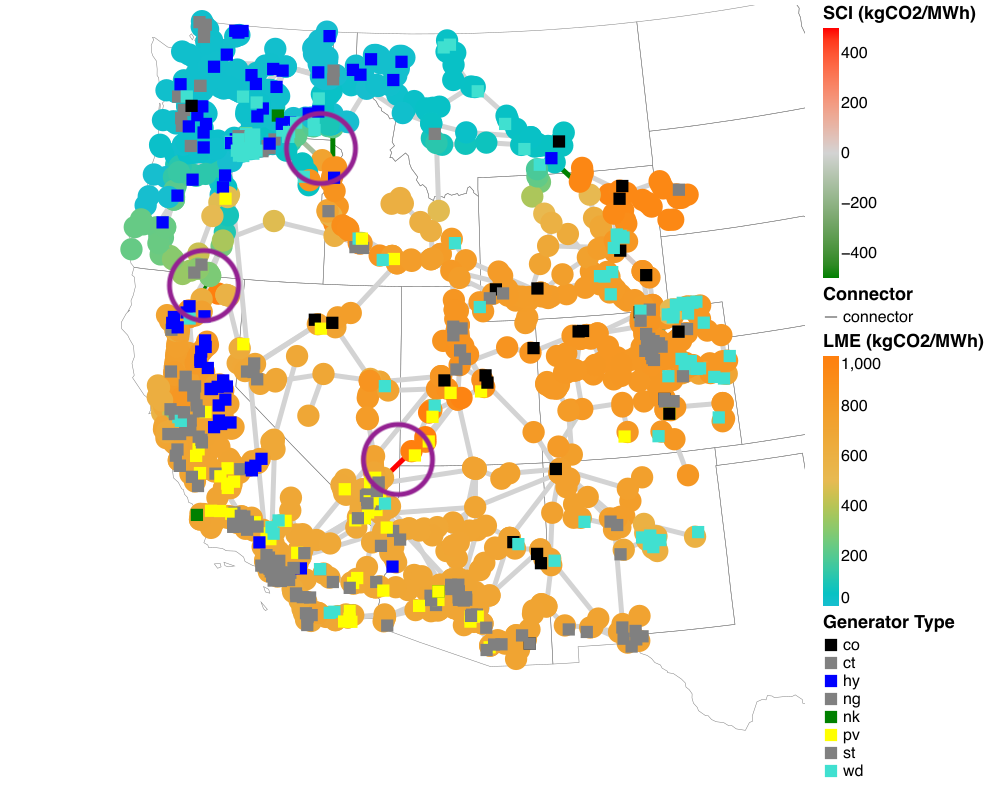}
    \label{fig:nightLME}
   \end{minipage}\hfill
   \begin{minipage}{0.48\textwidth}
     \centering
    \includegraphics[width=\linewidth]{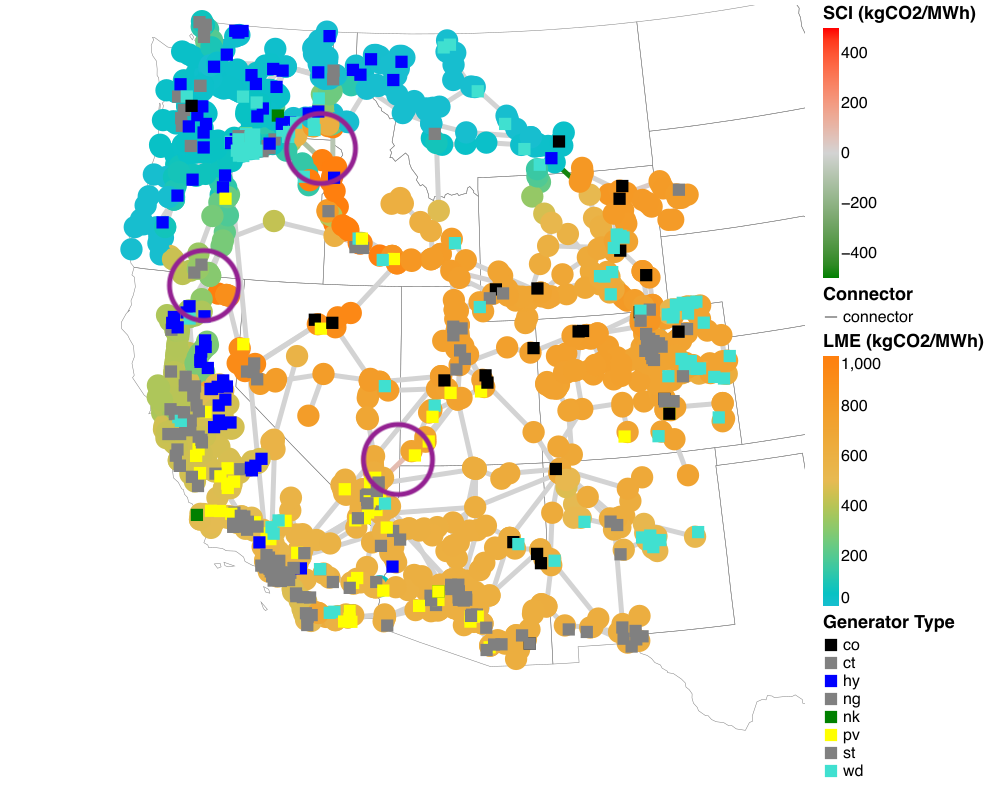}
    \label{fig:dayLME}
   \end{minipage}
   \caption{Average locational marginal emissions and shadow carbon intensities over the experimental year at 5am (left) and 2pm (right). Areas with transmission lines which experience significant hourly differences are circled in purple for clairty.}
    \label{fig:LME_time_dif}
\end{figure}

The hourly and seasonal trends of LMEs within the WECC grid remain well grouped into the three broad LME characteristic types. 

Due to the significantly higher temporal correlation of solar production when compared to wind, the hourly trends of LMEs across the WECC seem to be driven by the interplay of solar production and net load. In general, despite an overall increasing load through the afternoon as shown in \Cref{fig:hourly_demand}, the ``Inter-Mountain West'' and ``Sunbelt'' regions experience slight dips in average LMEs in \Cref{fig:hourly_LME}, likely due to the significant increase in solar production occurring at this time outpacing the demand increases. This increased solar production also leads to a decrease in net load in these regions, decreasing some of the congestion occurring in lines from the PNW region into northern California and into Idaho along with lines from Utah to southern Nevada as seen in \Cref{fig:LME_time_dif}.

\begin{figure}[h!]
   \begin{minipage}{0.48\textwidth}
     \centering
     \includegraphics[width=\linewidth]{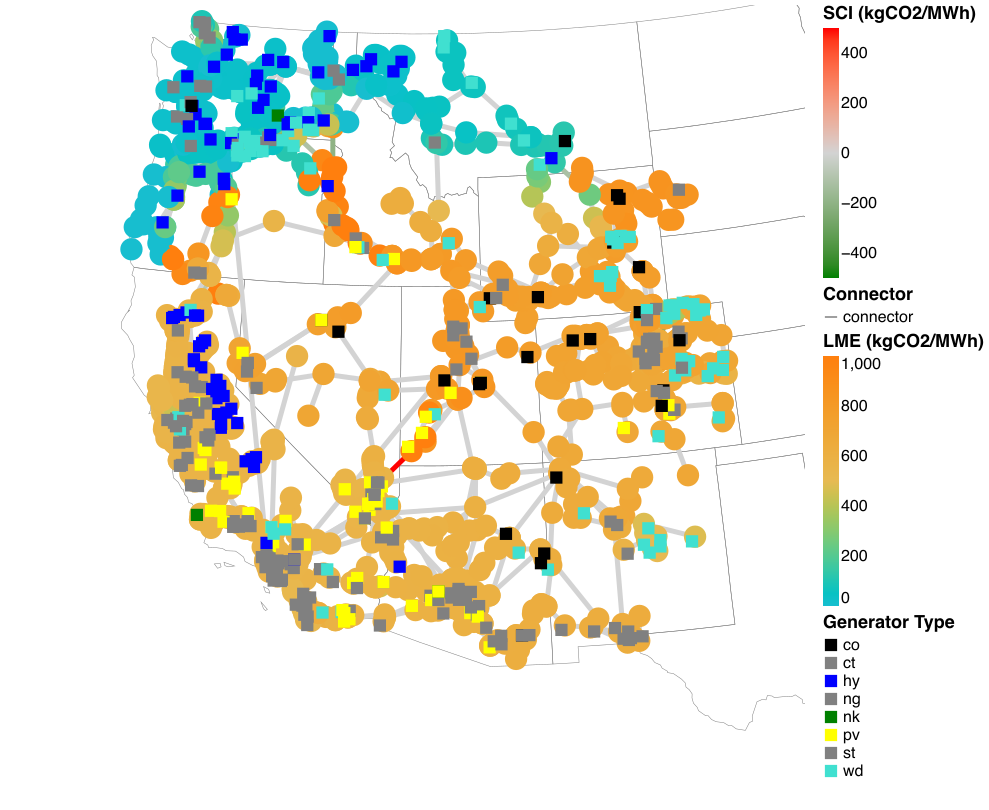}
    \label{fig:winterLME}
   \end{minipage}\hfill
   \begin{minipage}{0.48\textwidth}
     \centering
    \includegraphics[width=\linewidth]{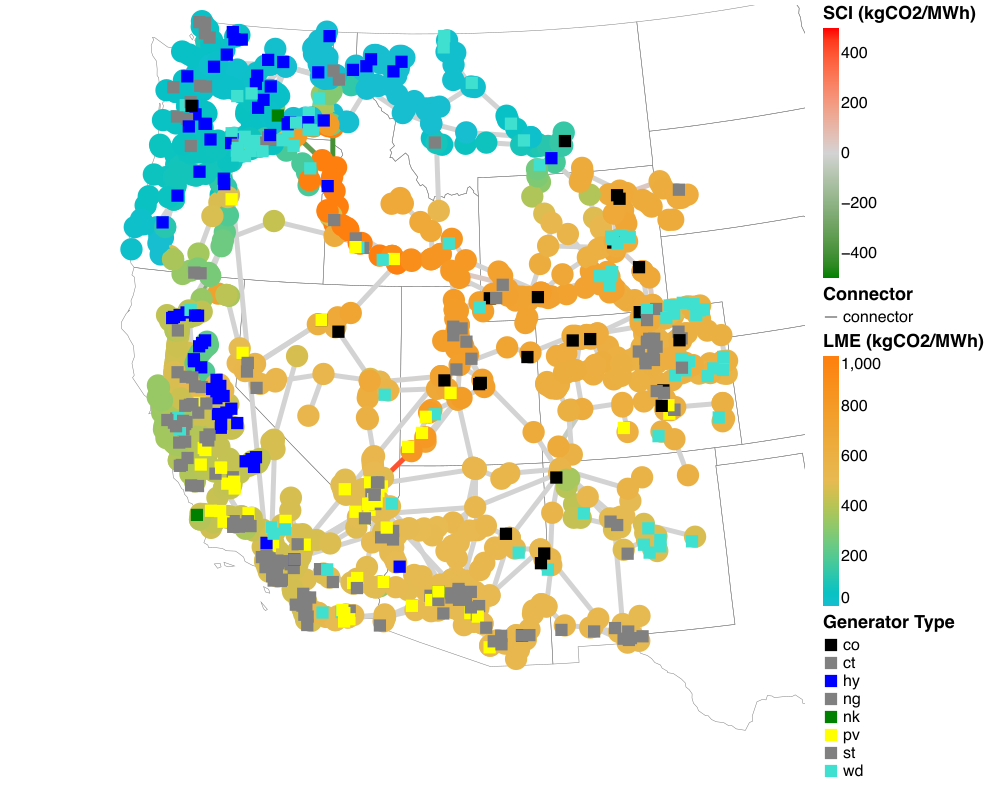}
    \label{fig:summerLME}
   \end{minipage}
   \caption{Average locational marginal emissions and shadow carbon intensities over the experimental month in January (left) and July (right).}
   \label{fig:LME_months_dif}
\end{figure}

The monthly LME trends are similarly dominated by solar production trends (see \Cref{fig:monthly_demand}-\Cref{fig:monthly_LME}), with higher solar production levels of the summer months leading to significantly decreased LMEs across the ``Inter-Mountain West'' and ``Sunbelt'' regions, despite a net increase in demand (likely driven by cooling systems) during these months. Interestingly, some southern portions of the PNW region also experience a slight decrease in LMEs during these months (see \Cref{fig:LME_months_dif}), likely driven by a combination of decreasing demand in the PNW region during the summer months and a decreased demand for energy exports from increased solar production in California. It is also worth noting that the Idaho region remains an exception to the trend of decreasing LMEs in the summer months, showing a regional average LME actually higher than that of coal plants, indicating that consumption in the region is likely causing increased fossil fuel dispatch and decreased renewable dispatch through a combination of regional congestion and the physics of power flows.


\subsection*{LMEs Balance System Emissions and Reflect Entity Level Impacts}\label{sec:accounting}
\begin{figure}[h!]
    \centering
    \includegraphics[width=0.7\linewidth]{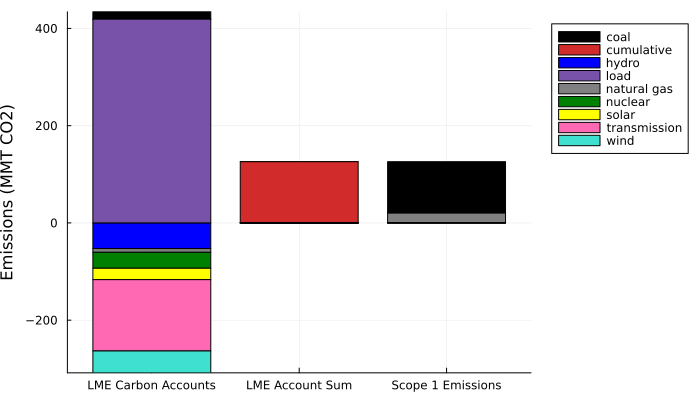}
    \caption{Total carbon accounts over the experimental horizon.}
    \label{fig:accounts}
\end{figure}
In this section, we apply the LME-based accounting scheme in the manner of \nameref{method:acct}. This scheme provides individual carbon accounts for all grid entities including loads, generation, and transmission. 

Over the course of the experiment, we observed a total of 125.9 million metric tons (MMT) of scope 1 carbon emissions from generation resources. Applying an LME-based accounting framework, we observe that this translates into a carbon account of 419.1 MMT ascribed to the grid load, a net total of -146.6 MMT ascribed to generation resources, and a total of -146.6 MMT ascribed to transmission infrastructure. This yields a total of 125.9 MMT of carbon emissions across all carbon accounts, directly matching the scope 1 emissions of the grid dispatch.
\begin{table}[h!]
    \centering
    \begin{tabular}{ |c||c|c|c|c|c|c| } 
     \hline
      Generation Type & Coal & Gas & Nuclear & Wind & Solar & Hydro\\ 
        \hline
      Carbon Account (MMT \chem{CO_2}) & 15.2 & -7.7 & -32.8 & -45.0 & -23.5 & -52.7\\ 
      Dispatch (GWh) & 122,131 &  41,876 & 64,666 & 86,153 & 39,744 & 268,392\\ 
      Scope 1 Emissions (MMT \chem{CO_2}) & 105.6 & 20.2 & - & - & - & -\\
     \hline
    \end{tabular}
    \caption{Generation resource carbon accounts and total power production over the experimental horizon.}
    \label{tab:gen_accounts}
\end{table}

\begin{figure}[h!]
   \begin{subfigure}{0.48\textwidth}
     \centering
     \includegraphics[width=\linewidth]{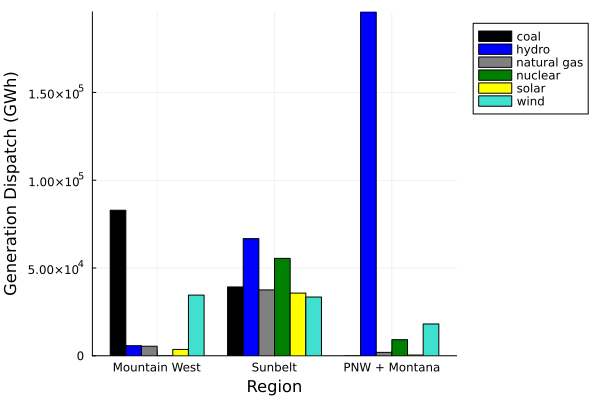}
     \caption{dispatch}
     \label{fig:region_dispatch}
   \end{subfigure}\hfill
   \begin{subfigure}{0.48\textwidth}
     \centering
    \includegraphics[width=\linewidth]{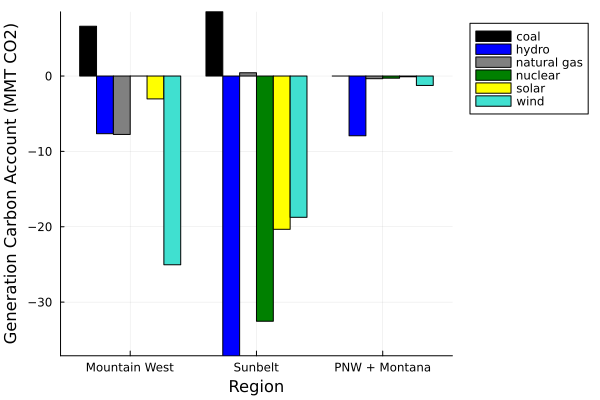}
    \caption{carbon accounts}
    \label{fig:region_accts}
   \end{subfigure}
   \caption{Total generation dispatch and carbon accounts by LME characteristic groups over the experiment year.} \label{fig:region_breakdown}
\end{figure}

The large magnitude of load-attributed LME-based emissions accounts is notable, with the total load carbon account summing to almost 3 times the scope 1 carbon emissions of the WECC grid. This is indicative of a grid with high amounts of renewable production without saturation, i.e., despite higher production levels, renewables do not spend much time acting as a marginal generation resource, causing the load to be assigned carbon emissions at a rate higher than the grid average emission factor\footnote{Grid average emission factor refers to the total emissions of the entire power grid system divided by the total power produced by the system.}. In this instance, we observe the average LME-based emissions account intensity to be $\approx$672 kg\chem{CO_2}/MWh, roughly halfway between the carbon intensities of natural gas and coal, indicating a general tendency of fossil fuel plants to be the marginal resources across the entirety of WECC, while the grid average emission factor is $\approx$202 kg\chem{CO_2}/MWh, or roughly halfway in between the carbon intensities of renewables and natural gas, indicating a significant presence of renewables in the overall generation mix.

Dividing the net generation account into each fuel type, we see that coal generation experiences the only positive carbon account with a total of 15.2 MMT carbon emissions while natural gas production received a negative carbon account of -7.7 MMT \chem{CO_2} despite being a scope 1 emissions source (\Cref{tab:gen_accounts}). However, in general, both coal and natural gas carbon accounts are quite small in magnitude compared to all other grid entities (\Cref{fig:accounts}), indicating that such fossil resources spend most of their time on the margin. 

Upon closer inspection of the generation accounting at the regional level we see that the negative carbon account of natural gas is driven by the coal-heavy ``Mountain West'' region, despite natural gas having a much larger dispatch footprint within the ``Sunbelt'' regions (\Cref{fig:region_breakdown}). Furthermore, we observe that the intensity of the mountain west natural gas carbon account is -1423 MT \chem{CO_2}/GWh, much higher than the difference between natural gas and coal generation carbon intensity ($\approx$350 MT \chem{CO_2}/GWh). When considering the dispatch cost ordering of generators within the simulation (coal being cheaper than gas), this indicates that the negative gas carbon accounts present in this region likely stems from natural gas plants being used to relieve regional congestion and thus enable an increased dispatch of renewables.

\begin{figure}[h!]
    \centering
    \includegraphics[width=0.6\linewidth]{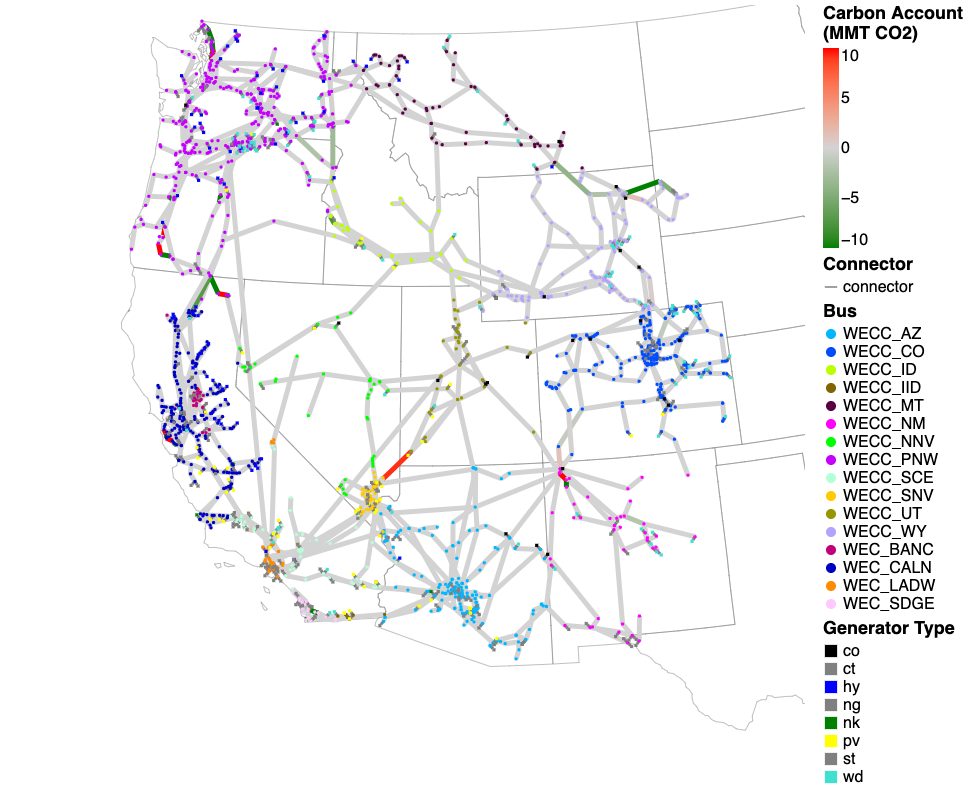}
    \caption{Total transmission carbon accounts over the experimental horizon.}
    \label{fig:line_accounts}
\end{figure}

Additionally, while, in total, the transmission accounts experience a significant net negative, we observe a distribution of positive and negative carbon accounts across transmission infrastructure. Lines whose congestion limits the dispatch of coal and necessitates an increased dispatch of natural gas receive positive carbon accounts. This indicates that their upgrade would allow for more coal generation, and thus lead to an increase in emissions. Lines whose congestion limits the dispatch of renewable resources receive negative carbon accounts, indicating that their upgrade would result in a decrease in emissions. Here, we once again observe lines connecting the low carbon PNW and Montana regions to California (darker green lines in the center-west of \Cref{fig:line_accounts}) and ``Intermountain West'' (lighter green lines connecting PNW to Idaho and Montana to Wyoming in \Cref{fig:line_accounts}) to experience negative carbon accounts, implying the existence of significant congestion and the possibility of carbon (and cost) reducing upgrades. Furthermore, we observe that congestion is not entirely renewable limiting and we observe some lines to experience positive carbon accounts - particularly those connecting the coal-margin ``Intermountain West'' to the often gas-margin ``Sunbelt'' regions (see red lines connecting Utah to Nevada and Colorado to New Mexico in \Cref{fig:line_accounts}). The existence of both positive and negative carbon accounts across the WECC transmission system suggests the importance of a principled upgrade strategy and the benefit that indicators such as the shadow carbon rent can provide for actors wishing to improve a transmission network to encourage decarbonization.

\subsection*{LMEs Provide Accurate Emission Signals for Data Center Operations}\label{sec:demand_response}
\paragraph{Carbon Signal Based Operation} In the following experiments, we evaluate the use of LMEs for carbon-based demand response programs. In particular, we consider a data center operator with eight data centers located in regions with high expected data center future load, as projected by NREL\cite{roberts_data_2025}. 
\begin{table}[h!]
    \small
    \centering
    \begin{tabular}{ |c||c|c|c|c|c|c|c|c| } 
     \hline
      \textbf{Location} & San Francisco & Los Angeles & Sacramento & Reno & Salt Lake & Phoenix & Denver & Seattle \\ 
      \hline
      \textbf{Bus \#} & 1167 & 67 & 800 & 292 & 25 & 205 & 322 & 28\\ 
      \textbf{Load \#} & 594 & 48 & 438 & 158 & 19 & 108 & 173 & 21\\ 
     \hline
    \end{tabular}
    \caption{Selected data center locations}
    \label{tab:datacenter_locations}
\end{table}
We perform four distinct experiments, varying the base load size of the data centers, the data center operating regions, and considering a grid-wide operator with control over all eight data centers vs. an operator localized to California, considering only the San Francisco, Los Angeles, and Sacramento data center locations. The explicit details of each experiment setup are summarized in \Cref{tab:shifting_results}. 

For all experiments, we assume that up to 20\% of a data center's load may be shifted as demand response. Shifted loads must be realized within the same day, and no data center load can exceed an additional 20\% of baseload at any time\footnote{The 20\% shifting assumption is seen in other academic work\cite{gorka_electricityemissionsjl_2025} and based on the results of data center load shifting trials at Oracle which reported an ability to decrease load by 25\%\cite{giacobone_google_2025}.}, i.e. for a 100MW baseload data center its load may be shifted to take on any value between 80MW and 120MW during any dispatch period. 

For each experiment the WECC dispatch is first reran with the addition of the experiment's data centers as static base loads, referred to as the ``base case''. The data center operator then solves a load shifting optimization problem, where they seek to minimize the sum of hourly LMEs at each data center node (given by the base case) multiplied by the data center loads at each hour. The difference between the hourly LMEs multiplied by the base load and the optimal objective of this problem is considered to be the expected shift in overall system emissions resulting from the application of such carbon-aware demand response
\begin{equation}
    \text{expected change}:= \sum_{t, dc} \text{load shift}_{t,dc} \times \text{LME}_{t,dc},
\end{equation}
while the total carbon emissions of the grid resulting from the economic dispatch of the shifted loads minus the realized carbon emissions from the economic dispatch of the base case describe the realized change in emissions
\begin{equation}
    \text{realized change}:= \text{shifted dispatch emissions} - \text{base case emissions}.
\end{equation}
\begin{table}
    \centering
    \begin{tabular}{|c||c|c|c|c|}
        \hline
        Experiment & Exp 1 & Exp 2 & Exp 3 & Exp 4 \\
        \hline
        Data Center Regions & WECC & CA & WECC & CA \\
        Data Centers & 8 & 3 & 8 & 3 \\
        Data Center Size (MW) & 100 & 100 & 200 & 200 \\
        \hline
        Expected Change (MT \chem{CO_2}) & -4,168 & -1,292 & -7,731 & -2,299\\
        Realized Change (MT \chem{CO_2}) & -3,833 & -1,151 & -6,860 & -1,971\\
        Change Ratio & 92\% & 89\% & 89\% & 86\%\\
        \hline
    \end{tabular}
    \caption{data center shifting experiment descriptions and results}
    \label{tab:shifting_results}
\end{table}
The results of each experiment are summarized in the bottom three rows in \Cref{tab:shifting_results}. In all experiments, LMEs provided an accurate estimate of the true system response with respect to data center load shifting, with the true response resulting in close to 90\% ratio of the LME-derived expected emissions change. 


Observing the differences between experiments, we see that the inclusion of geographically varied load shifting options has an outsized impact on the carbon-reduction ability of LME-based load-shedding schemes, with experiment 1 leading to a realized change of 333\% of that of experiment 2, despite having only 267\% more controllable load. This trend is repeated in experiments 3 and 4 with experiment 3 containing a realized shift of 348\% that of experiment 4. Intuitively, such advantages come from the additional arbitrage opportunity which exists between regions with differing LME characteristics, especially if such regions have differing temporal patterns. Of additional note is that the increase in data center size seems to have a diminishing impact, with 179\% and 171\% difference in realized change between experiments 1 and 3 and experiments 2 and 4, respectively, despite a doubling of the shiftable data center load. Here, such differences can likely be explained by the increase in system-wide LMEs caused by the increased overall load from the increased data center size leading to fewer time periods where renewable capacity is sufficient in comparison to load to act as the marginal resource.

\begin{figure}[h!]
   \begin{minipage}{0.48\textwidth}
     \centering
     \includegraphics[width=\linewidth]{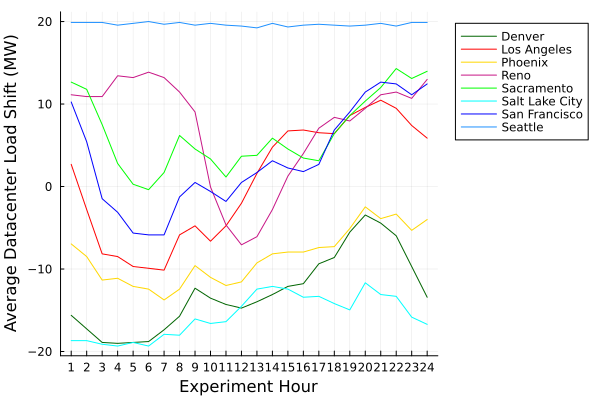}
   \end{minipage}\hfill
   \begin{minipage}{0.48\textwidth}
     \centering
    \includegraphics[width=\linewidth]{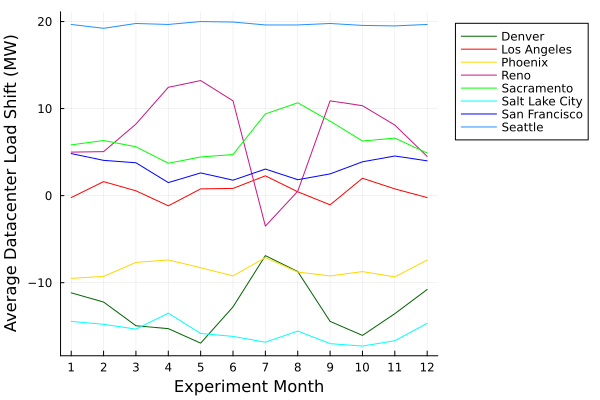}
   \end{minipage}
   \caption{Average hourly (left) and monthly (right) load shift factors for each data center in experiment 1.}
   \label{fig:dc_shifts}
\end{figure}
\begin{figure}[h!]
   \begin{minipage}{0.48\textwidth}
     \centering
     \includegraphics[width=\linewidth]{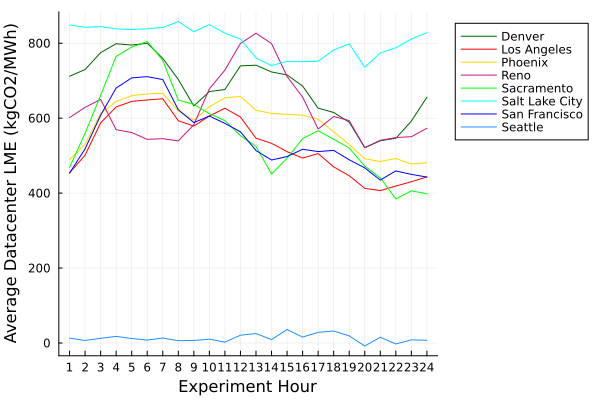}
   \end{minipage}\hfill
   \begin{minipage}{0.48\textwidth}
     \centering
    \includegraphics[width=\linewidth]{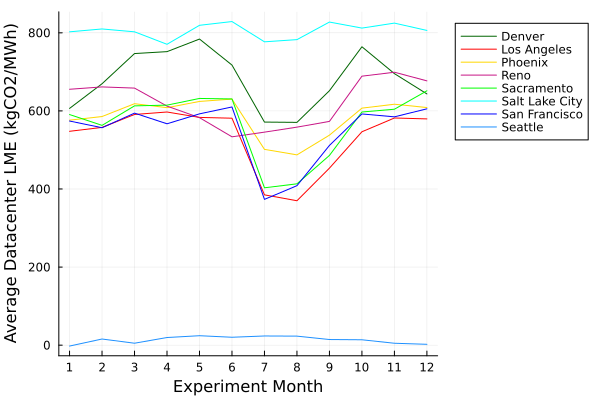}
   \end{minipage}
   \caption{Average hourly (left) and monthly (right) LMEs for each data center bus in experiment 1.}\label{fig:dc_LMEs}
\end{figure}

We can also consider the operational trends across the data centers within these experiments. As one might expect, the data center hourly load shifting patterns were highly correlated with the LME patterns for each region. The consistent low LMEs of the PNW led to a stable increased load for the Seattle data center over all experiment hours, while the rest of the data centers - located in the more variable ``Sun Belt'' and ``Intermountain-West'' regions shifted toward off-peak demand hours of the later night and early morning. Interestingly, we do not observe a significant shift towards off-peak midday hours as well, where high levels of solar production were shown to lower the LMEs of such regions. This is likely due to the urban, high-demand area siting of these data centers, causing them to be downstream of any congestion curtailing renewable production and thus less likely to have decreased LMEs during times of higher renewable production. Indeed, when observing the LME trends of the data center buses (\Cref{fig:dc_LMEs}) we observe that the pattern of midday LME drops from solar production does not exist for most buses with the exception of Sacramento. 

However, unlike the hourly case, when considering the monthly LME patterns of the data center buses we see the observed regional seasonal patterns remain true. Due to the 24-hour constrained shifting window, such seasonal patterns are not directly reflected in terms of load distribution. Rather, as the summer increase in solar production and subsequent decrease in LMEs has a differing magnitude of impact on different data center buses, the LME-based ordering of data centers changes and is thus reflected in the seasonal load shifting patterns.

\begin{figure}[h!]
   \begin{minipage}{0.48\textwidth}
     \centering
     \includegraphics[width=\linewidth]{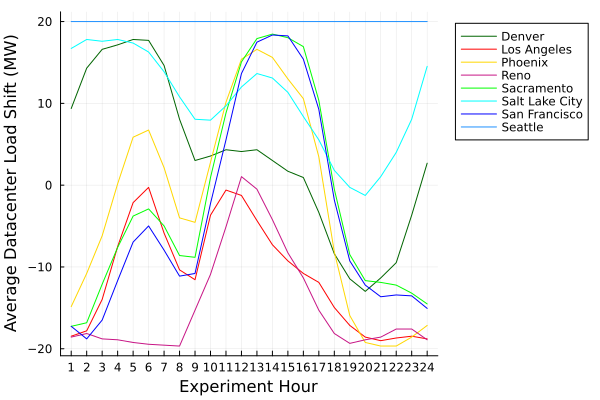}
   \end{minipage}\hfill
   \begin{minipage}{0.48\textwidth}
     \centering
    \includegraphics[width=\linewidth]{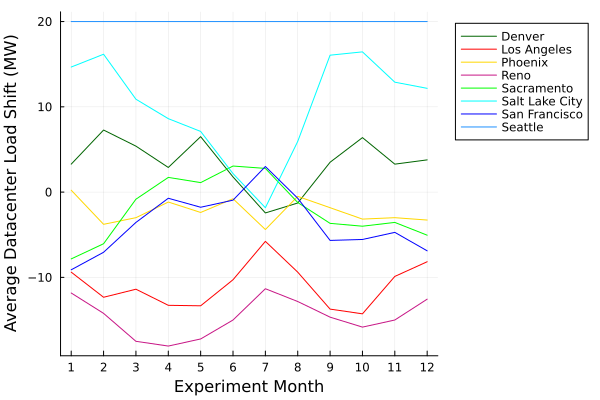}
   \end{minipage}
   \caption{Average hourly (left) and monthly (right) load shift factors for each data center when rerunning experiment 1 with locational marginal pricing based shifting.}
   \label{fig:dc_shifts_LMP}
\end{figure}

\paragraph{Economic Signal Based Operation } Data center operation is not constrained to only follow emissions signals, and is incentivized in many markets to follow economic signals given by locational marginal prices (LMPs). Previous literature has explored such data center operation and found that such operation may lead to higher system-level emissions. Indeed, when rerunning experiment 1 using LMPs as a shifting signal rather than LMEs we found that while system cost decreased by $\approx$5\%, system emissions increased by 78.5 MT \chem{CO_2}. These emission increases are likely driven by the economically incentivized shifting of power away from ``Sunbelt'' data centers and into ``Mountain West'' data centers due to the lower fuel price (yet higher emissions) of coal in comparison to natural gas which is often marginal in the ``Sunbelt'' regions. Such intuition is corroborated in \Cref{fig:dc_shifts_LMP} where we see significant shifts away from Los Angeles and Reno data centers and significant positive shifts into Salt Lake City and Denver data centers. 

\subsection*{LMEs Provide Accurate Signals for Data Center and Generation Siting} \label{sec:siting}

\paragraph{Data Center Siting} We consider the use of LMEs to understand the emissions impact of deciding siting static loads within a grid network. Here we randomly select a bus and a dispatch hour and observe how the overall system emissions change from adding a constant 200MW load, assumed to represent a large data center. The expected change in emissions at each experiment hour is given by the hourly LME multiplied by a 200MW load, while the realized change in system emissions is calculated by the difference between the overall system emissions with no additional load and the overall system emissions with the added 200MW load resulting from solving the DC optimal power flow dispatch. We perform this experiment for 1000 buses and dispatch periods chosen at random and observe the total realized change in emissions over the experimental horizon for all buses to be 1093.6 MT \chem{CO_2}, or 114\% of the expected change of 959.1 MT \chem{CO_2} from the total expected change in emissions, indicating that LMEs provide a fairly accurate estimate of system response for data center (and other similarly sized) load siting.

\paragraph{Generation Siting} Finally, we consider the use of LMEs to understand the emissions impact of deciding siting generation within a grid network. Here we repeat the methodology from static load siting; however rather than adding a 200MW load we add a negative 200MW load to a bus to simulate a dispatched generation resource sized in the fashion of a large wind or solar generation resource. We again performed this experiment for 1000 buses and dispatch periods chosen at random and observed the total realized change in emissions to be -845.0 MT \chem{CO_2}, totaling to 92.6\% of the total expected change in emissions of 912.8 MT \chem{CO_2}, again indicating LMEs provide an overall accurate estimate of system response for generation siting as well.

\begin{figure}[h!]
   \begin{minipage}{0.48\textwidth}
     \centering
     \includegraphics[width=\linewidth]{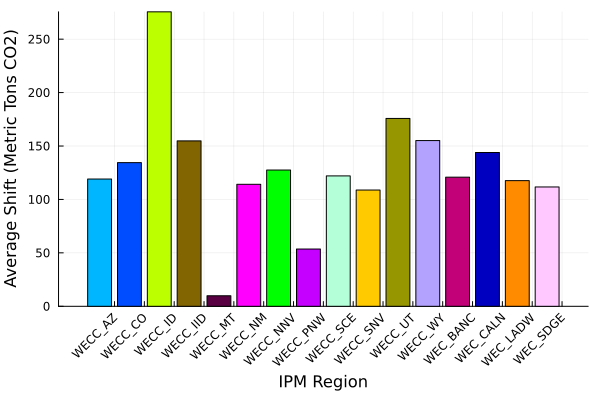}
    \label{fig:load_site_regions}
   \end{minipage}\hfill
   \begin{minipage}{0.48\textwidth}
     \centering
    \includegraphics[width=\linewidth]{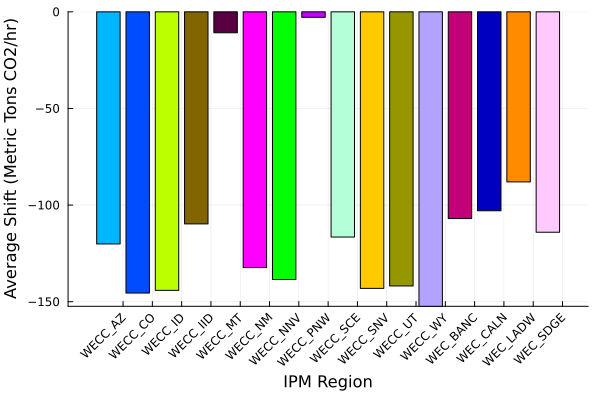}
    \label{fig:gen_site_regions}
   \end{minipage}
   \caption{Average emissions shift from 200MW load (left) and generation (right) siting within each IPM region.}
\end{figure}

The siting experiments also serve to confirm consumer insights of the LME patterns, with consumption sources sited in the PNW and Montana causing significantly lower additional emissions than those sited elsewhere. Additionally, the siting of generation sources in the ``Mountain West'' regions leads to some of the highest magnitudes of emission displacement. Indeed, extrapolating these results to MT \chem{CO_2}/MWh, a consumer building a 200 MW data center in the PNW would expect only to need to build a zero carbon resource that produces power at $\approx$ 67MW (assuming constant production over all hours) to have a net zero system level emission impact.

\section*{DISCUSSION}


\subsection*{Insights for Data Center and Renewable Developers}
The LME and accounting characteristics of the simulated WECC grid described in \nameref{sec:spatial}, \nameref{sec:temporal}, and \nameref{sec:accounting} provide significant insight into how data center operators can expect each of the three carbon intervention strategies to be most effective in reducing system emissions.

\begin{enumerate}
    \item \textit{The PNW and Montana regions have lowest LMEs for data center siting:} the stable, low LMEs of the PNW and Montana regions make such areas attractive options as operating sites for large energy consumers, such as data centers, looking to have a decreased carbon footprint. After these areas, the lower, sun-varied LMEs of ``Sunbelt'' regions, particularly those of northern California make these regions more attractive over the more consistent high-LME ``Inter-Mountain West'' regions.
    
    \item \textit{The Sunbelt regions are attractive for carbon-aware data center operation:} For data centers operating within the PNW and Montana regions, the stability of the low LMEs implies that there may be little to gain from consumption changes such as load shifting. However, for consumers operating outside these areas, particularly in the ``Sunbelt'' regions, the correlation between LME temporal patterns and solar strength indicates reliable opportunities to lower a data center's carbon footprint via carbon-aware operation. In particular, operational changes, which shift loads toward high solar production and subsequent low LME hours and months in these regions, have a high potential to provide significant carbon benefits.
    
    \item \textit{The Inter-Mountain West regions have strong potential for carbon reduction projects:} For developers looking to take carbon reduction actions, the ``Inter-Mountain West'' regions exhibit strong potentials for such projects. The high LMEs of these regions combined with low intra-region transmission congestion and high potential for VRE resources leads to likely strong system carbon reductions per unit of VRE capacity added within the regions. Transmission upgrades between the PNW and Montana regions and the ``Inter-Mountain West'' regions also exhibit strong potential for significant system carbon reductions through the highly negative SCIs exhibited across such lines.
\end{enumerate}

\subsection*{Insights for Policy Analysis}
\paragraph{LME-based interventions show promise.} There has been much policy discussion around scope 2 accounting systems to improve upon existing issues and better align carbon-conscious consumers such as the data center operators discussed above. Recent updates from the GHG protocol working group have signaled interest in the use of LMEs as an improved evaluation framework for consumer interventions\cite{huckins_scope_2025}.

Indeed, in Sections \nameref{sec:demand_response} and \nameref{sec:siting}, we observed that insight from LMEs proved accurate at providing system response estimates for load siting, demand response, and generation procurement actions, indicating that LME-based intervention strategies and LME-based accounting can be powerful tools for aligning data center consumption with the renewable transition.

\paragraph{Network representation and interconnection-scale simulation are important.} It should be noted that such accuracy findings contradict some previous studies \cite{gorka_electricityemissionsjl_2025}, where LME-based demand response was found to have significantly lower accuracy than observed here. A key difference between these studies was the size and scope of the network dispatch simulation. The larger scale WECC simulation in our study seems to provide more stability to regional LMEs with respect to significant load changes, comparing to simulations restricted to more local areas such as California. This implies that, to properly align intervention actions with system response, planners and operators must be careful to simulate their interventions on a power system with sufficiently large scale and high spatial and temporal resolution.

\paragraph{Fine-grained simulation provides insights on tax credit policy.} 

Furthermore, taking in conjunction the observed accuracy of LMEs to describe system response and the observed LME characteristics, we see the importance of fine-grained approaches to carbon accounting policies. LMEs showed significant variance across often considered blanket policy regions such as FERC Electricity Markets\cite{federal_energy_regulatory_commission_electric_2025} (for example, Idaho and the PNW are both within the same Northwest electricity market but experience drastically different LME characteristics), along with even individual nodes exhibiting significant temporal variation. Such policy implications are broadly applicable beyond just the WECC grid and can be used to evaluate the structure of current decarbonization incentives. One particularly pertinent example is the green hydrogen tax credit, where hydrogen producers operating in the ``qualifying states'' of Washington and California, having implemented greenhouse gas cap programs, can qualify for green hydrogen tax credits without additional energy procurement measures\cite{internal_revenue_service_credit_2025}. Such a broad-scope state-level policy fits poorly with the observed current LME trends where many nodes within California exhibit average LMEs close to the level of natural gas carbon intensity, indicating that, despite the existence of emissions caps, in the short term additional demand is met through increased fossil fuel consumption. Of course, such an analysis does not account for the full complexities of such a policy - in this instance likely including a desire to create political incentives for states to adopt emissions caps - however, it does provide important insight into how such policies may be misaligned with the physical grid response.

\subsection*{Recommendations for Regional Policy Makers.}
Finally, to complete our exploration, we can consider the insights from Sections \nameref{sec:spatial}, \nameref{sec:temporal}, and \nameref{sec:accounting} through the lens of regional policymakers under the adoption of such LME-based accounting systems.

\emph{Regions with high penetration of firm low carbon power:} For policymakers in these regions, the existing high frequency of low LMEs provides strong incentives for data centers to locate here and shift workloads to this area. Policy makers wishing to continue to provide such carbon incentives should ensure that the low carbon generation capacity remains sufficient to meet an increasing future load. 

\emph{Regions with consistent marginal fossil fuels:} Despite possibly containing a significant renewable capacity, the marginal responsibility of emissions from fossil fuel plants mean that data centers are disincentivized from locating workloads in these regions. Policy makers looking to attract data center operations should look to lower average LMEs through increasing the amount of time that carbon-free resources act as marginal generations through building out carbon-free generation capacity in the region. Such an increase could likely be achieved partially through carbon-free generation projects financed through data center development, as high regional LMEs likely make such areas attractive for energy procurement strategies.

\emph{Regions with high solar penetration:} The variable LMEs present in these regions provide an opportunity for policy makers to use data center intervention actions to flatten the infamous duck curve present in high solar regions. In particular, the lower LMEs observed during peak solar production incentivize data centers performing LME-based workload shifting interventions to shift demand into hours of peak solar production, creating an overall more stable load profile when considering net demand with VRE resources removed. Such interventions may help improve grid resilience while simultaneously allowing for increased renewable deployment.

\newpage


\section*{METHODS}\label{sec:methods}


\subsection*{Generation Dispatch}\label{method:dispatch}
For each experiment, we run 8760 separate power dispatch problems, each simulating one hour of WECC system operation. To do so, we follow the LME calculation methodology similar to that of \citet{rudkevich_locational_2012} and use a combined dispatch and LME calculation model for the WECC grid. In particular, we form a feasible set of possible generation dispatches from the set of optimal solutions to the DC-OPF economic dispatch problem with load shedding:
\begin{align}
   DCOPF(P^D):= \min_{P^G, \theta, f} \quad  & c_{gen}^T P^G + c_{shed} \cdot 
   \|ls\|_1 \tag{DCOPF Model} \label{p:dcOPF}\\
    \text{s.t.}\quad &
    B_{ij} (\theta_i - \theta_j) - f_\ell =  0, &\forall \ell=(i,j) \in L,\\
    & \sum_{g \in G_i} P^G_g - \sum_{\ell = (i,j) \in L}f_\ell + \sum_{\ell = (j,i) \in L}f_\ell - P_i^D+ls_i= 0, & \forall i \in N,\\
    & f_\ell \leq f^{max}_\ell, & \forall \ell \in L,  \\
    & f_\ell \geq -f^{max}_\ell, & \forall \ell \in L,  \\
    & P^G_g \leq P^{max}_g, & \forall g \in G,\\
    &P^G_g \geq P^{min}_g, & \forall g \in G, \\
    & ls_i \geq 0, & \forall i \in N, \\
    & ls_i \leq P^D_i, & \forall i \in N.
\end{align} 
We then optimize over the feasible set of optimal solutions to the \eqref{p:dcOPF} by minimizing the secondary objective of system carbon emissions, allowing for the calculation of marginal emission characteristics for each dispatch hour. The resulting optimal power flow solution from this two-tier optimization defines the dispatch solution for each hour.
\subsection*{LME-Based Carbon Accounting Calculation}\label{method:acct}
We apply the LME-based accounting scheme as described in \Cref{def:accts}. Such accounts are calculated matching the results of \citet{rudkevich_locational_2012} using the LME and SCI values derived from the generation dispatch for each hour along with the load, generation, and flow values associated with the corresponding dispatch. 

Additionally, for the most accurate picture of LME-based accounting within the WECC and to account for imperfections in our grid modeling methodology, we remove all high load shedding event hours from the accounting consideration, where high load shedding is considered to be an hour where load shedding is greater in magnitude than 100MWh.
\subsection*{Data Center Load Shifting}
For the load shifting experiments described in \nameref{sec:demand_response}, we run three distinct optimization problems. First, the considered data center base loads are added to the corresponding WECC nodes across all 8760 dispatch periods and we rerun \eqref{p:dcOPF} in the manner of \nameref{method:dispatch} to define emissions and LME values in the ``base case''. Second, we solve an emission minimizing load shifting problem across the data centers where each data center is allowed to increase or decrease its load by up to 20\% in any given dispatch period subject to the constraint that the sum of the total load over all participating data centers throughout each day (a 24-hour window) must remain the same as in the ``base case''. This leads to 365 separate data center load dispatch problems whose objective describes the total expected emissions change for each day:
\begin{align}
   CO2SHIFT(LME):= \min_{\Delta} \quad  & \sum_{t \in T, i \in DC}LME_{t,i} \Delta_{t,i} \tag{Shifting Model} \label{p:shift}\\
    \text{s.t.}\quad &
    \sum_{t \in T, i \in DC}\Delta_{t,i} = 0 \\
    & \Delta_{t,i} \leq 0.2 \times DCSIZE &\forall t \in T, i \in DC \\
    & -\Delta_{t,i} \leq 0.2 \times DCSIZE &\forall t \in T, i \in DC
\end{align} 
Finally, the resulting data center load shifts from the above problems are applied and the 8760 dispatch periods are reran with the shifted loads to define emissions in the ``shifted case''. We then subtract the total ``shifted case'' emissions from the ``base case'' emissions to generate the realized emissions change. 

\newpage


\section*{RESOURCE AVAILABILITY}


\subsection*{Lead contact}


Requests for further information and resources should be directed to and will be fulfilled by the lead contact, Andy Sun (sunx@mit.edu).

\subsection*{Materials availability}


This study did not generate new materials.

\subsection*{Data and code availability}

Any information required to reanalyze the data reported in this paper is available from the lead contact upon request. 

\section*{ACKNOWLEDGMENTS}


The authors would like to thank the gift from Meta that in part supported this work and for the valuable reviews and comments from Hank He and Nikky Avila of the Global Energy Team at Meta. Additionally, L.C. was partially supported by the Dick and Jerry Smallwood Fellowship Fund during portions of this work. The authors would also like to thank Thomas Lee for the compilation and provision of WECC grid network and demand data.

\section*{AUTHOR CONTRIBUTIONS}


Conceptualization, A.S. and L.C.; methodology, A.S. and L.C.; investigation, A.S. and L.C.; writing-–original draft, L.C.; writing-–review \& editing, A.S. and L.C.; funding acquisition, A.S.; resources, A.S. and L.C.; supervision, A.S.

\section*{DECLARATION OF INTERESTS}


The authors declare no competing interests.

\section*{DECLARATION OF GENERATIVE AI AND AI-ASSISTED TECHNOLOGIES}


During the preparation of this work, the authors used Google Gemini in order to generate code scaffolding for plots and figures. After using this tool or service, the authors reviewed and edited the content as needed and take full responsibility for the content of the publication.






\newpage

\bibliography{references}

\bigskip


\newpage

\end{document}